\documentclass[12pt]{amsart}

\usepackage{pstcol}
\psset{unit=0.5cm,linewidth=0.4pt}
\def\vertex{\pscircle[fillstyle=solid,fillcolor=black]{0.1}}
\definecolor{light}{gray}{0.9}
\definecolor{medium}{gray}{0.8}

\let\Bbb=\mathbb
\let\frak=\mathfrak
\def\RR{{\Bbb R}}
\def\QQ{{\Bbb Q}}
\def\ZZ{{\Bbb Z}}

\def\mm{{\frak m}}
\def\nn{{\frak n}}

\def\CC{{\mathcal C}}
\def\DD{{\mathcal D}}
\def\fpn{{\fp n}}
\def\fp#1{{^{[#1]}}}
\def\ce#1{{\lceil #1\rceil}}
\def\iso{\cong}
\def\fon{{\fo n}}
\def\fo#1{{^{1/#1}}}

\let\tensor=\otimes
\def\ol#1{\overline{#1}}

\let\oldhash=\#
\def\#{\mathbin{\oldhash}}

\def\chara{\operatorname{char}}
\def\Hom{\operatorname{Hom}}

\def\depth{\operatorname{depth}}
\def\vol{\operatorname{vol}}
\def\gp{\operatorname{gp}}
\def\Cl{\operatorname{Cl}}
\def\HK{\textup{HK}}

\newtheorem{lemma}{Lemma}[section]
\newtheorem{corollary}[lemma]{Corollary}
\newtheorem{theorem}[lemma]{Theorem}
\newtheorem{proposition}[lemma]{Proposition}

\theoremstyle{definition}

\newtheorem{remark}[lemma]{Remark}

\begin{document}
\title{Conic divisor classes over a normal monoid algebra}
\author{Winfried Bruns}
\address{Universit\"at Osnabr\"uck,
FB Mathematik/Informatik, 49069 Osnabr\"uck, Germany}
\email{Winfried.Bruns@mathematik.uni-osnabrueck.de}

\maketitle

In \cite{BGu:Div} J. Gubeladze and the author have studied the
divisorial ideals of an algebra $R=K[M]$ where $K$ is a field and
$M$ a normal affine monoid. The divisorial ideals which have a
monomial $K$-basis cut out from the group $\gp(M)$ by a translate
of the cone $\RR_+M$ have been called \emph{conic}. Up to
isomorphism the conic ideals are exactly the direct summands of
the extension $R^\fon$ of $R$ where $R^\fon$ is the $K$-algebra
over the monoid $(1/n)M=\{x/n: x\in M\}\subset \QQ\tensor\gp(M)$.

In this note we want to extend the discussion in \cite{BGu:Div}.
In a rather straightforward manner we will represent the divisor
classes by the full-dimensional open cells in a certain cell
complex. The conic classes can then be identified with the
full-dimensional open cells in a decomposition of the torus
$\bigl(\RR\tensor\gp(M)\bigr)/\gp(M)$. Furthermore, these classes
can be characterized by the relative compactness of a certain
group associated with them. For torsion classes this group is
finite.

By Hochster's theorem $R^\fon$ is a Cohen-Macaulay ring. Therefore
conic classes are Cohen-Macaulay. Using a criterion of St\"uckrad
and Vogel \cite{StV} for the Cohen-Macaulay property of Segre
products, Bae\c tica \cite{Ba} has given examples of
Cohen-Macaulay divisorial ideals that are not conic. We will
review his construction and streamline the arguments somewhat.

In the last part of the paper we investigate the multiplicities of
the conic classes in the decomposition of $R^\fon$, as $n$ varies.
This multiplicity turns out to be a quasi-polynomial for all $n\ge
1$, given by the number of lattice points in the union of the
interiors of certain cells of the cell complex mentioned above. As
Watanabe \cite{Wa} has already noticed, this argument can be used
for the computation of the Hilbert-Kunz multiplicity of $R$ in
characteristic $p>0$. In addition it yields some assertions about
the Hilbert-Kunz function of $R$.

\section{Conic divisor classes}
Let $M$ be a normal affine monoid of rank $d$. The divisor class
group of $K[M]$ has been discussed in \cite{BGu:Div}. We refer the
reader to this article for details that may be missing in the
following. We always assume that $M$ is positive, i.~e.\ $0$ is
the only invertible element in $M$.

Changing the embedding if necessary, we may assume that
$\ZZ^d=\gp(M)$. One has $M=\ZZ^d\cap \RR_+M$ since $M$ is normal,
and the cone $\RR_+M$ has an irredundant representation
$$
\RR_+M=\{x\in\RR^d: \sigma_i(x)\ge 0,\ i=1,\dots,s\}
$$
as an intersection of halfspaces defined by primitive integral
forms $\sigma_i$, $i=1,\dots,s$ on $\RR^d$. (``Primitive
integral'' means that the elements of the canonical basis of
$\RR^d$ have coprime integral values.)

Every (fractional) divisorial ideal of $R$ is isomorphic to a
(fractional) monomial ideal of $R$, and the latter are precisely
the ideals
$$
\DD(u)=K\cdot\{x\in\ZZ^d: \sigma(x)\ge u\}\qquad u\in\RR^s.
$$
Here $\sigma(x)=(\sigma_1(x),\dots,\sigma_s(x))$, and the
inequality must be read componentwise. (We identify monomials with
their exponent vectors.) Of course, $\DD(u)=\DD(\ce u)$ with $\ce
u=(\ce{u_1},\dots, \ce{u_s})$.

The ideals $\DD(u)$ and $\DD(v)$ are isomorphic as $R$-modules if
and only if there exists a monomial $x\in\ZZ^d$ such that
$\DD(u)=x\DD(v)$. For $u,v\in\ZZ^s$ this amounts to
$u=v+\sigma(x)$ for some $x\in\ZZ^d$. Therefore the divisor class
group $\Cl(R)$ is given by $\ZZ^s/\sigma(\ZZ^d)$.

Some more information is contained in the topological group
$\RR^s/\sigma(\ZZ^d)$ and its cell decomposition $\bar\Delta$
defined as follows. First, let $\Delta$ be the cell decomposition
of $\RR^n$ defined by all the hyperplanes through integral points
that are parallel to the coordinate hyperplanes. Then the open
$s$-cells $\delta$ of $\Delta$ are in $1$-$1$-correspondence to
$\ZZ^s$ via their ``upper right'' corners
$$
\ce{x}\in\ZZ^s,\qquad x\in\delta.
$$
´The ``upper closures''
$$
\ce{\delta}=\{y\in\RR^s: \ce y=\ce x\},\qquad x\in\delta,
$$
are the fibers of the map $\RR^s\to\ZZ^s$, $y\mapsto \ce y$.
Clearly $\ce{\delta}$ is a semi-open cube and the union of open
cells of $\Delta$. Moreover, these semi-open cubes correspond
bijectively to the divisorial ideals of $M$. The cell
decomposition is invariant under the operation of $\sigma(\ZZ^d)$
on $\RR^s$ via translations, and so $\Delta$ induces a cell
decomposition $\bar\Delta$ of $\RR^s/\sigma(\ZZ^d)$. Collecting
all our arguments, we obtain

\begin{theorem}
The following sets are in a natural bijective correspondence:
\begin{itemize}
\item[(a)] $\Cl(R)$;
\item[(b)] the fibers of the map $\RR^s/\sigma(\ZZ^d)\mapsto
\Cl(R)$ induced by the assignment $x\mapsto \ce x$, $x\in\RR^s$;
\item[(c)] the open $s$-cells of $\bar\Delta$.
\end{itemize}
\end{theorem}

As pointed out in the introduction, the \emph{conic} divisorial
ideals
$$
\CC(y)=K\cdot \bigl(\ZZ^d\cap (y+\RR_+M)\bigr), \qquad y\in\RR^d.
$$
are of special interest. Clearly,
$$
\CC(y)=\DD(\sigma(y))
$$
so that the conic divisorial ideals are exactly those among the
$\DD(u)$ for which $u$ can be chosen in $\sigma(\RR^d)$.

If one restricts the cell decomposition of $\RR^s$ to $\RR^d\iso
\sigma(\RR^d)$, then one obtains the cell decomposition $\Gamma$
of $\RR^d$ by the hyperplanes
$$
H_{i,z}=\{x\in\RR^d:\sigma_i(x)=z\},\qquad i=1,\dots,s,\ z\in\ZZ.
$$
It induces a cell decomposition $\bar\Gamma$ of the torus
$\RR^d/\ZZ^d$, which we can also define by restricting
$\bar\Delta$ to the submanifold $\RR^d/\ZZ^d$ (via $\sigma$).
Figure \ref{CellDec13} gives an example of such a decomposition.
We have marked the cone $\RR_+M$ and the semiopen square
$(-1,0]^2$, a fundamental domain for the action of $\ZZ^2$.
Evidently $\bar\Gamma$ has exactly $3$ open cells of dimension
$2$.
\begin{figure}[htb]
\begin{center}
\psset{unit=1cm}
\def\vertex{\pscircle[fillstyle=solid,fillcolor=black]{0.05}}
\begin{pspicture}(-1,-1)(2,3)
 \pspolygon[fillstyle=solid, fillcolor=light,linewidth=0pt](0,0)(1,3)(2,3)(2,0)
 \pspolygon[fillstyle=solid, fillcolor=light,linewidth=0pt](0,0)(-1,0)(-1,-1)(0,-1)
 \multirput(-1,-1)(1,0){4}{\multirput(0,0)(0,1){5}{\vertex}}
 \multirput(-1,-1)(0.3333,0){6}{\psline(0,0)(1.333,4)}
 \multirput(-1,-1)(0,1){5}{\psline(0,0)(3,0)}
 \rput(-1,0){\psline(0,0)(1,3)}
 \rput(-1,1){\psline(0,0)(0.666,2)}
 \rput(-1,2){\psline(0,0)(0.333,1)}
 \rput(1,-1){\psline(0,0)(1,3)}
 \rput(1.333,-1){\psline(0,0)(0.666,2)}
 \rput(1.666,-1){\psline(0,0)(0.333,1)}
 \psline[linewidth=1.0pt](0,0)(2,0)
 \psline[linewidth=1.0pt](0,0)(1,3)
\end{pspicture}
\end{center}
\caption{} \label{CellDec13}
\end{figure}

\begin{corollary}\label{tor-dec}
The following sets are in a natural bijective correspondence:
\begin{itemize}
\item[(a)] the set of conic divisor classes;
\item[(b)] the fibers of the map $\RR^d/\ZZ^d\to
\RR^s/\sigma(\ZZ^d)$ induced by the assignment $x\mapsto \lceil
\sigma(x)\rceil$;
\item[(c)] the full dimensional cells of $\bar\Gamma$.
\end{itemize}
\end{corollary}

Note that $\RR_+M$ contains a point $x$ with $\sigma_i(x)>0$ for
all $i$. This implies that the full-dimensional open cells of
$\Gamma$ are contained in full-dimensional open cells of $\Delta$.

We can now characterize the conic classes in terms of the sizes of
certain subgroups of $\RR^s/\sigma(\ZZ^d)$:

\begin{corollary} A divisor class is
\begin{itemize}
\item[(a)] a torsion element if and only if it
is represented by an ideal $\DD(u)$ for which $\ZZ \bar u$ is a
finite subgroup of $\RR^s/\sigma(\RR^d)$;
\item[(b)] conic if it
is represented by an ideal $\DD(u)$ for which $\ZZ \bar u$ is a
relatively compact subgroup of $\RR^s/\sigma(\ZZ^d)$ (i.~e.\ the closure
of $\ZZ\bar u$ is compact).
\end{itemize}
\end{corollary}

The next proposition connects the property of being conic with a
condition of Stanley \cite[p.\ 41]{Sta}.

\begin{proposition}\label{char-conic}
The divisorial ideal $\DD(u)$ is conic if and only if the coset
$\sigma(\RR^d)-\ce u$ in $\RR^s/\sigma(\RR^d)$ contains a point in
the semi-open cube $Q=(-1,0]^s$.
\end{proposition}

\begin{proof}
We can replace $u$ by $\ce u\in\ZZ^s$, and therefore assume that
$u\in\ZZ^s$.

Suppose first that $\DD(u)=\CC(x)$. Then $u=\ce{\sigma(x)}$ and
$\sigma(x)-u=\sigma(x)-\ce{\sigma(x)}$ is in
$Q\cap(\sigma(\RR^d)-u)$.

Conversely, suppose that $v\in(Q\cap\sigma(\RR^d))-u$. Then
$u+v\in \sigma(\RR^d)$, and $\DD(u)=\DD(u+v)$.
\end{proof}

Let $T$ be the torsion subgroup of $\Cl(R)\iso
\ZZ^s/\sigma(\ZZ^d)$. Then $\Cl(R)/T$ is a free abelian group of
rank $m=s-d$. It can be identified with a rank $m$ lattice $L$ in
the $m$-dimensional vector space $\RR^s/\sigma(\RR^d)$. Roughly
speaking, the conic classes form the set of lattice points in a
polytope in this vector space.

\begin{corollary}\label{con-cl-con}
Let $P$ be the polytope in $\RR^s/\sigma(\RR^d)$ spanned by the
images of the conic classes in $\Cl(R)/T$. Then
\begin{itemize}
\item[(a)] a class $c\in\Cl(R)$ is conic if and only if its
residue class modulo $T$ belongs to $P$;
\item[(b)] $P$ has dimension $m=s-d$, and is centrally symmetric
with respect to $[\omega]/2$ modulo $T$ where $\omega$ is the
canonical module of $R$.
\end{itemize}
\end{corollary}

\begin{proof}
Let $P'$ be the image of the semi-open cube $Q$ in
$\RR^s/\sigma(\RR^d)$. By the proposition a class is conic if and
only if its image modulo $T$ belongs to $P'$. Therefore the convex
hull $P$ of the images of the conic classes has property (a). It
is of dimension $m$ since it contains the origin of $L$ (given by
the trivial class) and the classes of a system of generators of
$\Cl(R)$: the classes of the monomial divisorial prime ideals
generate $\Cl(R)/T$ and are conic (see \cite[3.4]{BGu:Div}).
Therefore $P$ must have dimension $m$. The central symmetry with
respect to $[\omega]/2$ results from the fact that $[\omega]-c$ is
conic if $c$ is conic; see loc.\ cit.
\end{proof}

\section{Conic and Cohen-Macaulay classes}

As pointed out above, the conic classes represent Cohen-Macaulay
modules. In general however, there exist Cohen-Macaulay divisorial
monomial ideals that are not conic, as shown by Bae\c tica
\cite{Ba}, using a theorem of St\"uckrad and Vogel on the
Cohen-Macaulay property of Segre products. Bae\c tica applies it
to divisorial ideals over Segre products of $3$ polynomial rings.
However, he does not use the full strength of the results in
\cite{StV}, and his arguments can be improved.

Let $R_i$, $i=1,2$, be a positively graded, finitely generated
$K$-algebra with irrelevant maximal ideal $\mm_i$. Furthermore,
let $M_i$ be a graded module over $R_i$, $i=1,2$ such that $\depth
M_i\ge 2$. Then the $k$-th local cohomology of the module $M=M_1\#
M_2$ over the Segre product $R=R_1\#R_2$ (with respect to the
irrelevant maximal ideal $\mm$ of $R$) is the direct sum
\begin{equation}
\bigoplus_{i+j=k+1} \bigl(H_{\mm_1}^i(M_1)\#H_{\mm_2}^j(M_2)\bigr)
\oplus \bigl(M_1\#H_{\mm_2}^k(M_2)\bigr)\oplus
\bigl(H_{\mm_1}^k(M_1)\#M_2\bigr).\tag{$*$}
\end{equation}
This is just the K\"unneth formula for Serre cohomology
\cite[0.2.10]{StV}, rewritten in terms of local cohomology. (See
\cite[p.\ 38]{StV} for the correspondence between Serre and local
cohomology.) Since $\depth M_i\ge 2$, one has $H_{\mm_i}^k(M_i)=0$
for $k=0,1$, and so $\depth M\ge 2$ by \cite[0.2.12]{StV}.

In order to control the vanishing of the local cohomology we
introduce the following notation: for a graded maximal
Cohen-Macaulay module $N$ over a positively graded affine
$K$-algebra $S$ of dimension $e$, with irrelevant maximal ideal
$\nn$, we set
\begin{align*}
b(N)&=\inf\{k:N_k\neq 0\},\\
h(N)&=\sup\{k:H_\nn^e(N)_k\neq 0\}.
\end{align*}
Under certain conditions $N$ and $H_{\nn}^{e}(N)$ are ``gap
free'': $N_k\neq0$ for all $k\ge b(N)$ and $H_\nn^e(N)_k\neq 0$ or
all $k\le h(N)$. By graded local duality, $H_{\nn}^{e}(N)$ is the
graded $K$-dual of the maximal Cohen-Macaulay module
$\Hom_S(N,\omega)$ where $\omega$ denotes the canonical module of
$S$. Therefore it is enough to discuss gap freeness for $N$. If
the graded component $S_1$ is not contained in a minimal prime
ideal of $S$, then it is not contained in an associated prime
ideal of $N$, and after the extension of $K$ to an infinite field
it follows that $S_1$ contains a nonzerodivisor of $N$, whence $N$
is gap free.

It is now easy to generalize \cite[Th.\ I.4.6]{StV} to more than
two factors in the Segre product:

\begin{proposition}\label{Segre-CM}
Suppose that $R_1,\dots,R_n$ are positively graded Cohen-Macaulay
$K$-algebras of dimensions $d_1,\dots,d_n\ge 2$, with irrelevant
maximal ideals $\mm_1,\dots,\mm_n$, and let $M_i$ be a graded
maximal Cohen-Macaulay module over $R_i$, $i=1,\dots,n$. If
$(R_i)_1$ is not contained in a minimal prime ideal of $R_i$ for
$i=1,\dots,n$, then the following are equivalent:
\begin{itemize}
\item[(a)] $M_1\#\cdots\#M_n$ is Cohen-Macaulay;
\item[(b)] for all nonempty, proper subsets $I\subset\{1,\dots,n\}$ one
has
$$
\min\{h(M_i):i\in I\}<\max\{b(M_j):j\notin I\};
$$
\item[(c)] there exists a permutation $j_1,\dots,j_n$ of
$1,\dots,n$ such that $M_{j_1}\#\cdots\#M_{j_t}$ is Cohen-Macaulay
over $R_{j_1}\#\cdots\#R_{j_t}$, $t=1,\dots,n$.
\end{itemize}
\end{proposition}

\begin{proof}
An iterative use of the formula $(*)$ allows one to compute the
local cohomology of the Segre product. The $k$-th local cohomology
of $M=M_1\#\cdots\#M_n$ over $R=R_1\#\cdots\#R_n$ (with respect to
the irrelevant maximal ideal $\mm$ of $R$) is the direct sum of
the $R$-modules
$$
H_{\mm_{i_1}}^{d_{i_1}}(M_{i_1})\#\cdots
\#H_{\mm_{i_t}}^{d_{i_t}}(M_{i_t})\# M_{i_{t+1}}\#\cdots\# M_{i_n}
$$
where $0\le t\le n$, $1\le i_1<\dots <i_t\le n$, $i_{t+1}<\dots
<i_n$, $\{i_1,\dots,i_s\}= \{1,\dots,s\}$, and
$$
k=d_{i_1}+\dots+d_{i_t}-(t-1).
$$
In addition to $(*)$ this computation only uses that
$H_{\mm_i}^j(M_i)=0$ for $j<d_i$.

The equivalence of (a) and (b) follows immediately, and only the
implication (b)$\implies$(c) needs to be proved. We have to find a
permutation $j_1,\dots,j_n$ such that
\begin{align*}
h(M_{j_1}\#\cdots\#M_{j_t})&=\min(h(M_{j_i}: 1\le i\le t)
<b(M_{j_{t+1}})\\
b(M_{j_1}\#\cdots\#M_{j_t})&=\max(b(M_{j_i}): 1\le i\le
t)>h(M_{j_{t+1}})
\end{align*} for $t=1,\dots,n-1$.

Reordering $1,\dots,n$ if necessary, we can assume that
$b(M_1)\le\dots\le b(M_n)$. We set $u_1=1$. Then we split the
sequence $2,\dots,n$ into blocks
$$
B_v=\{u_v+1,\dots,u_{v+1}\}, \qquad v=1,\dots,w,
$$
where $h(M_i)\ge b(M_{u_v})$, $i=u_v+1,\dots,u_{v+1}-1$, but
$h(M_{u_{v+1}})<b(M_{u_v})$. Such a decomposition exists by (b),
applied successively to the subsets $\{1,\dots,u_v\}$. Each of the
blocks $B_1,\dots,B_w$ is then cyclically permuted to
$$
B_v'=\{u_{v+1},u_v+1,\dots,u_{v+1}-1\}
$$
and the desired permutation is finally given by the concatenation
$1,B_1',\dots,B_w'$. It is not hard to check that this permutation
indeed satisfies the desired inequalities.
\end{proof}

Let $R_1,\dots,R_n$ be polynomial rings over $K$ of dimensions
$d_1,\dots,d_n\ge 2$ with the standard grading by total degree.
Then their Segre product $S=R_1\#\cdots \# R_n$ is a monoid ring
over the Segre product
$$
\ZZ_+^{d_1}\#\cdots\#\ZZ_+^{d_n}
$$
where the Segre product construction is simply applied to the
monoids of monomials. It is not hard to see that
$\Cl(S)\iso\ZZ^{n-1}$, and that the divisorial ideals over $S$ are
represented by the Segre products of shifted copies of the $R_i$,
$$
D=R_1(-s_1)\#\cdots \# R_n(-s_n), \qquad (s_1,\dots,s_n)\in\ZZ^n.
$$
The divisor class of $D$ is completely determined by the
differences $(s_2-s_1,\dots,s_n-s_1)\in\ZZ^{n-1}$. One has
$b(R_i(-s_i))=s_i$ and $h(R_i(-s_i))=s_i-d_i$.

Note that the proof of the implication (b)$\implies$(c) contains
an algorithm by which one can check the Cohen-Macaulay property of
the Segre product. It is especially simple if all the dimensions
$d_i$ are equal of constant value $d$: if $s_1\le\dots\le s_n$,
then $D$ is Cohen-Macaulay if and only if $s_{i+1}<s_i+d$ for
$i=1,\dots,n-1$, and the other cases are obtained by permutation.
For $n=3$ and $d_1=d_2=d_3=3$ the Cohen-Macaulay classes are
indicated in Figure \ref{CM333}. The conic classes correspond to
the lattice points in the polytope given by the inequalities
$-2\le x,y,y-x\le 2$. In particular there exist more
Cohen-Macaulay classes than conic ones.
\begin{figure}[htb]
\begin{center}
\begin{pspicture}(-4,-4)(4,4)
 \pspolygon[fillstyle=solid, fillcolor=light,linecolor=white,
 linewidth=0pt](0,2)(0,0)(2,0)(4.5,2.5)(4.5,4.5)(2.5,4.5)
 \pspolygon[fillstyle=solid, fillcolor=light,linecolor=white, linewidth=0pt](0,2)(0,-2)(-4.5,-2)(-4.5,2)
 \pspolygon[fillstyle=solid, fillcolor=light,linecolor=white, linewidth=0pt](-2,0)(2,0)(2,-4.5)(-2,-4.5)
 \pspolygon[fillstyle=solid, fillcolor=medium, linewidth=0.8pt](-2,0)(0,2)(2,2)(2,0)(0,-2)(-2,-2)
 \multirput(-2,-4)(1,0){1}{\vertex}
 \multirput(-2,-3)(1,0){2}{\vertex}
 \multirput(-4,-2)(1,0){7}{\vertex}
 \multirput(-3,-1)(1,0){6}{\vertex}
 \multirput(-2,0)(1,0){5}{\vertex}
 \multirput(-2,1)(1,0){6}{\vertex}
 \multirput(-2,2)(1,0){7}{\vertex}
 \multirput(1,3)(1,0){2}{\vertex}
 \multirput(2,4)(1,0){1}{\vertex}
 \psline{->}(0,-4.5)(0,4.5)
 \psline{->}(-4.5,0)(4.5,0)
\end{pspicture}
\end{center}
\caption{} \label{CM333}
\end{figure}

It follows immediately from the computation of the local
cohomology that in the case $d_1=\dots=d_n=d$ only the depths
$kd-(k-1)$ are possible for $k=1,\dots,n$. For $n=3$, $d=3$, the
(unbounded) areas of the classes of depth $2d-1$ have been
shadowed.

\begin{remark}
Bae\c tica \cite{Ba} has also investigated Segre products
$T=R^{(c)}\#S^{(d)}$ of Veronese subalgebras of polynomial rings
$R$ and $S$. If $\gcd(c,d)=1$, then $\Cl(T)\iso \ZZ$. For example,
if $\dim R=\dim S=2$ and $c=3$, $d=2$, then the set of
Cohen-Macaulay classes is $\{-5,-3,\dots,5,7\}$ (with respect to a
suitable choice of generator) and does not form an interval in
$\Cl(T)$.
\end{remark}

\section{The decomposition of $R\fon$}

Let $M$ be an affine monoid and $G=\gp(M)$. The monoid $(1/n)M$
contains $M$, and so $R\fon$ contains $R$ as a subalgebra.
Moreover, $\gp((1/n)M)=(1/n)G$. For each residue class $c\in
(1/n)G/G$ we set
$$
I_c=R\fon\cap K\cdot c.
$$
Then $I_c$ is an $R$-submodule of $R\fon$, and $R\fon$ decomposes
into the direct sum of its $R$-submodules $I_c$, $c\in (1/n)G/G$.

Now suppose that $M$ is normal. Then $I_c=\RR_+M\cap(M+c)$, and
the parallel translation by $-c$ yields
$$
I_c\iso \CC(-c).
$$
Therefore $I_c$ is a conic divisorial ideal.

We consider the cell complex $\bar\Gamma$ on the torus
$\RR^d/\ZZ^d$ (after the identification of $\gp(M)$ and $\ZZ^d$).
For each open $d$-cell $\bar\gamma$ of $\bar\Gamma$ we choose an
open $d$-cell $\gamma$ of $\Gamma$ representing it and set
$$
\ce{\gamma}=\{x\in \RR^d:\ce{\sigma(x)}=\sigma(y)\},\qquad
y\in\gamma.
$$
Then $\ce{\gamma}$ is the union of finitely many open cells of
$\Gamma$ contained in the closure of $\gamma$. Moreover, all the
conic ideals $\CC(x)$, $x\in\ce{\gamma}$, coincide. In the
following we write $\CC_\gamma$ for $\CC(x)$, $x\in\ce{\gamma}$.

\begin{theorem}\label{con-mult}
Let $v_\gamma(n)$ be the multiplicity with which the isomorphism
class of $\CC_\gamma$ occurs in the decomposition of $R\fon$ as an
$R$-module. Then
$$
v_\gamma(n)=\#\biggl(\ce{\gamma}\cap \frac 1n\ZZ^d\biggr).
$$
In particular, there exists a quasi-polynomial
$q_\gamma:\ZZ\to\ZZ$ with rational coefficients such that
$v_\gamma(n)=q_\gamma(n)$ for all $n\ge 1$. One has
$$
q_\gamma(n)=\vol(\gamma)n^d+a_\gamma n^{d-1}+\tilde
q_\gamma(n),\qquad n\in\ZZ,
$$
where $a_\gamma\in\QQ$ is constant and $\tilde q_\gamma$ is a
quasi-polynomial of degree $\le d-2$.
\end{theorem}

\begin{proof}
The direct summand $I_c$ of $R\fon$, $c\in (1/n)G/G$, is
isomorphic to $\CC_\gamma$ if and only if
$(-c+\ZZ^d)\cap\ce{\gamma}\neq\emptyset$. Moreover,
$(-c+\ZZ^d)\cap\ce{\gamma}$ contains at most one point, since
$\ce{\sigma(x)}\neq\ce{\sigma(y)}$ if $x-y\in \ZZ^d$, $x\neq y$.
This proves the first assertion.

For the second we decompose $\ce{\gamma}$ into the open cells
$\sigma$ of $\Gamma$ of which it is the union. The generating
function
$$
H_\sigma(t)=\sum_{n=0}^\infty \#\biggl(\sigma\cap\frac
1n\ZZ^d\biggr) t^n
$$
is a rational function of degree $0$. In fact, it is the Hilbert
series of the canonical module of the Ehrhart ring of the rational
polytope $\overline\sigma$, and therefore has degree $0$. Summing
over the cells $\sigma$ we obtain that
$$
H_\gamma(t)=\sum_{n=0}^\infty v_\gamma(n) t^n
$$
is a rational function of degree (at most) $0$. Therefore
$v_\gamma(n)$ is given by a quasi-polynomial $q_\gamma$ for all
$n\ge 1$.

Only the $d$-dimensional cell $\gamma$ itself contributes to the
degree $d$ term of $q_\gamma$. By Ehrhart reciprocity the number
of points in $\gamma\cap (1/n)\ZZ$  is $(-1)^d\bar q(-n)$ where
$\bar q$ is the Ehrhart quasi-polynomial of the closure
$\ol{\gamma}$ of $\gamma$. Since $\ol{\gamma}$ is a
full-dimensional rational polytope, the leading coefficient of
$\bar q_\gamma$ is constant and equal to the volume of
$\ol{\gamma}$.

Next we have to show that the coefficient of $n^{d-1}$ in
$q_\gamma$ is constant. It gets contributions from $\bar \gamma$
and from the Ehrhart quasi-polynomials of the open $(d-1)$-cells.
The coefficient of $n^{d-1}$ in $\bar q_\gamma$ is constant since
the affine hulls of the facets of $\ol{\gamma}$ contain lattice
points (Stanley \cite[p.\ 237]{Sta}). In fact, the support forms
of $\RR_+M$ map $\ZZ^d$ surjectively onto $\ZZ$. For the same
reason the coefficient of $n^{d-1}$ is constant in the Ehrhart
polynomials of the $(d-1)$-cells since these are the interiors of
rational polytopes in the affine hulls of the facets of
$\ol{\gamma}$ (which, as stated already, contain lattice points).
\end{proof}

In this proof we have used results of Ehrhart and Stanley on
lattice points in rational polytopes and Hilbert functions; for
example, see \cite{Sta} or 4.4.5, 6.3.5 and 6.3.11 in \cite{BH}.

As a corollary we obtain a description of the Hilbert-Kunz
function of $K[M]$ for normal $M$ and a nice formula for the
Hilbert-Kunz multiplicity of a graded algebra with normalization
$K[M]$. The argument has similarly been used by Watanabe
\cite{Wa}. Also see Conca \cite{Co} and Eto \cite{Eto} for results
on the Hilbert-Kunz multiplicity of monoid rings. By $\mu_R(N)$ we
denote the minimal number of generators of an $R$-module $N$.

By definition, the \emph{Hilbert-Kunz multiplicity} of a finitely
generated graded module $B$ over a positively graded algebra $A$
over a field $K$ of prime characteristic $p$ is
$$
e_{\HK,A}(B)= \lim_{e\to\infty} \frac{\dim_K
B/\mm\fp{p^e}B}{p^{de}}
$$
where $\mm$ is the irrelevant maximal ideal of $A$ generated by
all elements of positive degree.

\begin{corollary}\label{HK}
Let $K$ be a an algebraically closed field of characteristic $p>0$
and $M$ a normal affine monoid of rank $d$.

\begin{itemize}
\item[(a)] Set $R=K[M]$ and let $\mm$ be the maximal ideal of $R$
generated by the monomials different from $1$. Then
$$
\dim_K R/\mm\fp{p^e}=\mu_R(R\fo{p^e})=\sum_{\gamma}
\mu_R(\CC_\gamma) v_\gamma(p^e),\qquad e\in\ZZ,\ e\ge 1,
$$
is the value of the quasi-polynomial $q_\HK=\sum_\gamma
\mu_R(\CC_\gamma)v_\gamma$ at $p^e$. It has constant leading
coefficient
$$
e_{\HK, R}(R)=\sum_\gamma \mu_R(\CC_\gamma)\vol(\gamma)\in\QQ,
$$
and also the coefficient of its degree $d-1$ term is constant and
rational.

\item[(b)] Suppose $S$ is a graded subalgebra of $R$ such that
$R$ is a finite $S$-module and the normalization of $S$. Then
$e_{\HK,S}(S)=\sum_{\gamma} \vol(\gamma)\mu_S (\CC_\gamma)\in\QQ$.
\end{itemize}
\end{corollary}

\begin{proof}
One has the isomorphism $\bar R/\mm\fp{p^e}\iso \bar R\fo{p^e}/\mm
\bar R\fo{p^e}$ of $K$-vector spaces, and so Nakayama's lemma
yields
$$
\dim_K R/\mm\fp{p^e}=\mu_R(\bar R\fo{p^e})=\sum_{\gamma}
v_\gamma(p^e)\mu_R(\CC_\gamma).
$$
Now all the assertions in (a) follow immediately from the theorem.

For (b) we note the exact sequence
$$
0\to S\to R\to T\to 0
$$
in which $T$ has Krull dimension $<d$. Therefore
$e_{\HK,S}(S)=e_{\HK,S}(R)$. The decomposition of $R\fo{p^e}$ over
$R$ is, a fortiori, a decomposition over $S$, and the formula for
$e_{\HK, S}(R)$ follows again from the theorem.
\end{proof}

That the coefficient of the degree $d-1$ term is constant, as we
have derived from combinatorial arguments, is in fact true in much
more generality; see Huneke, McDermott, and Monsky \cite{HMM}.

Corollary \ref{HK} allows one to give non-trivial examples of
algebras $K[M]$ for which the quasi-polynomial $q_\HK$ is a true
polynomial. The Ehrhart quasi-polynomial of a polytope with
integral vertices is a polynomial. Therefore it is sufficient that
all cells of $\Gamma$ have integral vertices, and this holds if
the support forms $\sigma_1,\dots,\sigma_s$ of the cone $C=\RR_+M$
form a totally unimodular configuration: every linearly independent subset of
$d$ elements is a basis of $(\ZZ^d)^*$.

\end{document}